\documentclass{amsart}[11pt]
\DeclareMathSymbol{\twoheadrightarrow}
{\mathrel}{AMSa}{"10}

\def\Q{{\mathbf Q}}
\def\Z{{\mathbf Z}}
\def\C{{\mathbf C}}
\def\R{{\mathbf R}}
\def\F{{\mathbf F}}
\def\f{{\mathfrak f}}
\def\R{{\mathfrak R}}

\def\CC{{\mathfrak C}}
\def\cc{{\mathfrak c}}

\def\Sn{{\mathbf S}_n}
\def\An{{\mathbf A}_n}

\def\Gal{\mathrm{Gal}}
\def\Perm{\mathrm{Perm}}

\def\Tr{\mathrm{Tr}}

\def\End{\mathrm{End}}
\def\Aut{\mathrm{Aut}}
\def\Hom{\mathrm{Hom}}

\def\cl{\mathrm{cl}}

\def\Lie{\mathrm{Lie}}
\def\I{\mathrm{Id}}

\def\GL{\mathrm{GL}}

\def\dim{\mathrm{dim}}

\def\P{{\mathbf P}}
\def\n{{n}}

\newtheorem{thm}{Theorem}[section]
\newtheorem{lem}[thm]{Lemma}

\newtheorem{prop}[thm]{Proposition}
\theoremstyle{definition}

\newtheorem{ex}[thm]{Example}
\newtheorem{exs}[thm]{Examples}

\newtheorem{rem}[thm]{Remark}

\title[Cyclic covers, jacobians and endomorphism rings]
{The endomorphism rings of jacobians of cyclic covers of the
projective line}
\author[Yuri G. Zarhin]{By Yuri G. Zarhin $\dagger$}
\thanks{ $\dagger$ Partially supported by the NSF}
\begin{document}
\maketitle

 \centerline{\sl Department of Mathematics, Pennsylvania
State University,}

 \centerline{\sl  University Park, PA 16802,USA}

\centerline{{\sl e-mail:} {\tt zarhin\char`\@math.psu.edu}}

 \begin{abstract}
Suppose $K$ is a field of characteristic zero, $K_a$ is its algebraic closure, $f(x) \in K[x]$ is an irreducible polynomial of degree $n \ge 5$, whose Galois group coincides either with the full symmetric group $\Sn$ or with the alternating group $\An$. Let $p$ be an odd prime, $\Z[\zeta_p]$ the ring of integers in the $p$th cyclotomic field $\Q(\zeta_p)$. Suppose $C$ is the smooth projective model of the affine curve $y^p=f(x)$ and $J(C)$ is the jacobian of $C$. We prove that the ring $\End(J(C))$ of $K_a$-endomorphisms of $J(C)$ is canonically isomorphic to $\Z[\zeta_p]$.
 \end{abstract}


\section{Introduction}
We write $\Z,\Q,\C$  for the ring of integers, the field of
rational numbers and the field of complex numbers respectively.
Recall that a number field is called a CM-field if it is a purely
imaginary quadratic extension of a totally real field.
 Let $p$ be an odd prime, $\zeta_p \in \C$ a primitive $p$th
root of unity, $\Q(\zeta_p)\subset \C$ the $p$th cyclotomic field
and $\Z[\zeta_p]$ the ring of integers in $\Q(\zeta_p)$. It is
well-known that $\Q(\zeta_p)$ is a CM-field of degree $p-1$.
We write $\F_p$ for the finite field consisting of $p$ elements.

Let $f(x) \in \C[x]$ be a polynomial of degree $n\ge 4$ without multiple roots.
 Let $C_{f,p}$ be a smooth projective model of the smooth affine curve
$$y^p=f(x).$$
It is well-known that the genus $g(C_{f,p})$  of $C_{f,p}$ is
$(p-1)(n-1)/2$ if $p$ does not divide $n$ and
 $(p-1)(n-2)/2$ if it does.
The map
 $$(x,y) \mapsto (x, \zeta_p y)$$
gives rise to a non-trivial birational automorphism
$$\delta_p: C_{f,p} \to C_{f,p}$$
of period $p$.

The jacobian $J^{(f,p)}:= J(C_{f,p})$ of $C_{f,p}$ is an abelian
variety of dimension  $g(C_{f,p})$. We write $\End(J^{(f,p)})$ for
the ring of endomorphisms of $J^{(f,p)}$ over  $\C$. By Albanese
functoriality, $\delta_p$ induces an automorphism of $J^{(f,p)}$
which we still denote by $\delta_p$; it is known
(\cite[p.~149]{Poonen}, \cite[p.~458]{SPoonen}) that
$$\delta_p^{p-1}+\cdots +\delta_p+1=0$$
in $\End(J^{(f,p)})$. This gives us an embedding
$$\Z[\zeta_p] \cong \Z[\delta_p] \subset \End(J^{(f,p)})$$
(\cite[p.~149]{Poonen}, \cite[p.~458]{SPoonen}).

Our main result is the following statement.

\begin{thm}
\label{endo}
Let $K$ be a subfield of $\C$ such that all the coefficients of $f(x)$
lie in $K$. Assume also that $f(x)$ is an irreducible polynomial in $K[x]$
of degree $n \ge 5$
 and its Galois group over $K$ is either the symmetric group $\Sn$ or
the alternating group $\An$. Then
 $$\End(J^{(f,p)})=\Z[\delta_p]
\cong\Z[\zeta_p].$$
In particular, $J^{(f,p)}$ is a simple
complex abelian variety.
\end{thm}

\begin{rem}
In the case when $p$ is a Fermat prime the assertion of Theorem
\ref{endo} is proven in \cite{ZarhinCrelle}. (Also in
\cite{ZarhinCrelle} the author  proved  that if the conditions of
Theorem \ref{endo} hold true then $\Z[\delta_p]$ is a maximal
commutative subring in $\End(J^{(f,p)})$ for all odd primes $p$.
See \cite{ZarhinSb} for a similar result in positive
characteristic when $p\mid n$ and $n\ge 9$.) The ``opposite" case
when $J^{(f,p)}$ is an abelian variety of CM-type was studied in
\cite{JN}. An analogue of Theorem \ref{endo} for hyperelliptic
jacobians (i.e., the case of $p=2$) was proven in \cite{Zarhin}
(see also \cite{ZarhinM}, \cite{ZarhinMRL2}).
\end{rem}

\begin{exs}
\begin{enumerate}
\item
 the polynomial $x^n-x-1 \in \Q[x]$ has Galois group $\Sn$ over
$\Q$ (\cite[p.~42]{Serre}). Therefore the endomorphism ring (over
$\C$) of the jacobian $J(C)$ of the  curve $C:y^p=x^n-x-1$ is
$\Z[\zeta_p]$ if  $n\ge 5$.
\item
the  Galois group of the ``truncated exponential"
 $$\exp_n(x):=1+x+\frac{x^2}{2}+\frac{x^3}{6}+ \cdots +
\frac{x^n}{n!}\in \Q[x]$$ is either $\Sn$ or $\An$ \cite{Schur}.
Therefore the endomorphism ring (over $\C$) of the jacobian $J(C)$
of the  curve $C:y^p=\exp_n(x)$ is $\Z[\zeta_p]$ if  $n\ge 5$.
\end{enumerate}
\end{exs}

\begin{rem}
\label{ribcom}
 If $f(x) \in K[x]$ then the curve $C_{f,p}$ and
its jacobian $J^{(f,p)}$ are defined over $K$. Let $K_a\subset
\C$ be the algebraic closure of $K$. Clearly, all endomorphisms
of $J^{(f,p)}$ are defined over $K_a$.
 This implies that in order to prove Theorem \ref{endo}, it suffices to
check that $\Z[\delta_p]$ coincides with the ring of all
$K_a$-endomorphisms of $J^{(f,p)}$ or equivalently, that
$\Q[\delta_p]$ coincides with the $\Q$-algebra of
$K_a$-endomorphisms of $J^{(f,p)}$.
\end{rem}

The paper is organized as follows.  Section \ref{MT} contains
auxiliary results about endomorphism algebras of complex abelian
varieties. We use them in Section \ref{Prelim} in order to study
endomorphisms of $J^{(f,p)}$.  In  Section \ref{prf} we prove the
main result. The short last Section contains corrigendum to
\cite{ZarhinCrelle}.

The author would like to thank the referee for useful comments.

\section{Complex abelian varieties}
\label{MT}
Throughout this section we assume that $Z$ is a complex abelian variety of positive dimension.
As usual, we write $\End^0(Z)$ for the semisimple
finite-dimensional $\Q$-algebra $\End(Z)\otimes \Q$. We write $\CC_Z$ for the center of $\End^0(Z)$.
 It is well-known that $\CC_Z$ is  a direct product of finitely many number fields.
  All the fields involved are either totally real number fields or CM-fields.
  Let $H_1(Z,\Q)$ be the first rational homology group of $Z$; it is a $2\dim(Z)$-dimensional $\Q$-vector space.
  By functoriality $\End^0(Z)$ acts on $H_1(Z,\Q)$; hence we have an embedding
$$\End^0(Z) \hookrightarrow \End_{\Q}(H_1(Z,\Q))$$
(which sends $1$ to $1$).

Suppose $E$ is a subfield of $\End^0(Z)$ that contains the
identity map. Then $H_1(Z,\Q)$ becomes an $E$-vector space of
dimension $$d=\frac{2\dim(Z)}{[E:\Q]}.$$ We write
$$\Tr_E: \End_E(H_1(Z,\Q)) \to E$$ for the corresponding trace
map on the $E$-algebra of $E$-linear operators in $H_1(Z,\Q)$.

 Extending by
$\C$-linearity the action of $\End^0(Z)$ and of $E$ on the
complex cohomology group $$H_1(Z,\Q)\otimes_{\Q}\C=H_1(Z,\C)$$ of
$Z$ we get the embeddings $$E\otimes_{\Q}\C \subset
\End^0(Z)\otimes_{\Q}\C \hookrightarrow \End_{\C}(H_1(Z,\C))$$
which provide $H_1(Z,\C)$ with a natural structure of free
$E_{\C}:=E\otimes_{\Q}\C$-module of rank $d$.
 If $\Sigma_E$ is the set of embeddings of $\sigma: E \hookrightarrow \C$
 then it is well-known that
$$E_{\C}=E\otimes_{\Q}\C=\prod_{\sigma\in \Sigma_E}
E\otimes_{E,\sigma}\C=\prod_{\sigma\in \Sigma_E}\C_{\sigma}$$
where
$$\C_{\sigma}=E\otimes_{E,\sigma}\C=\C.$$
Since $H_1(Z,\C)$ is a free $E_{\C}$-module of rank $d$, there is
the corresponding trace map
$$\Tr_{E_{\C}}:\End_{E_{\C}}(H_1(Z,\C)) \to E_{\C}$$
which coincides on $E_{\C}$ with multiplication by $d$ and with
$\Tr_E$ on $\End_E(H_1(Z,\Q))$.

We write $\Lie(Z)$ for the tangent space of $Z$; it is a
$\dim(Z)$-dimensional $\C$-vector space. By functoriality,
$\End^0(Z)$ and therefore $E$ acts on $\Lie(Z)$. This provides
$\Lie(Z)$ with a natural structure of $E\otimes_{\Q}\C$-module. We
have
$$\Lie(Z)=\bigoplus_{\sigma\in
\Sigma_E}\C_{\sigma}\Lie(Z)=\oplus_{\sigma\in
\Sigma_E}\Lie(Z)_{\sigma}$$ where
$$\Lie(Z)_{\sigma}=\C_{\sigma}\Lie(Z)=\{x \in \Lie(Z)\mid
ex=\sigma(e)x \quad \forall e\in E\}.$$ Let us put
$$n_{\sigma}=n_{\sigma}(Z,E)=\dim_{\C_{\sigma}}\Lie(Z)_{\sigma}=\dim_{\C}\Lie(Z)_{\sigma}.$$
We write $\sigma': E \hookrightarrow \C$ for the composition of
$\sigma: E \hookrightarrow \C$ and the complex conjugation $\C
\to \C$. The embedding $\sigma'\in \Sigma_E$ is (called) the
complex-conjugate of $\sigma$.

\begin{rem}
\label{conjugate}
 It is well-known (\cite[p.~53]{Deligne}, \cite[p.~84]{MZ}) that
 $$n_{\sigma}+n_{\sigma'}=d \quad \forall \sigma \in \Sigma_E.$$
 \end{rem}

\begin{rem}
\label{dual}
 Let $\Omega^1(Z)$ be the space of the differentials
of the first kind on $Z$. It is well-known that the natural map
$$\Omega^1(Z) \to \Hom_{\C}(\Lie(Z),\C)$$ is an isomorphism. This
isomorphism allows us to define via duality the natural
homomorphism
$$E  \to
\End_{\C}(\Hom_{\C}(\Lie(Z),\C))=\End_{\C}(\Omega^1(Z)).$$ This
provides $\Omega^1(Z)$ with a natural structure of
$E\otimes_{\Q}\C$-module in such a way that
$$\Omega^1(Z)_{\sigma}:=\C_{\sigma}\Omega^1(Z)\cong
\Hom_{\C}(\Lie(Z)_{\sigma},\C).$$ In particular,
$$n_{\sigma}=\dim_{\C}(\Lie(Z)_{\sigma})=\dim_{\C}(\Omega^1(Z)_{\sigma}).$$
\end{rem}

\begin{thm}
\label{mult} Suppose $E$ contains $\CC_Z$. Then the tuple
$$(n_{\sigma})_{\sigma\in \Sigma_E} \in \prod_{\sigma\in
\Sigma_E}\C_{\sigma}=E\otimes_{\Q}\C$$ lies in
$\CC_Z\otimes_{\Q}\C$. In particular, if $E/\Q$ is Galois and
$\CC_{Z} \ne E$ then there exists a nontrivial automorphism
$\kappa: E \to E$ such that $n_{\sigma}=n_{\sigma\kappa}$ for all
$\sigma\in \Sigma_E$.
\end{thm}

{\sl Proof.} The inclusion $\CC_Z \subset E$ implies that $\CC_Z$
is a field.

 There is a canonical Hodge decomposition (\cite[chapter 1]{MumfordAV}, \cite[pp.~52--53]{Deligne})
$$H_1(Z,\C)=H^{-1,0} \oplus H^{0,-1}$$
where $H^{-1,0}$ and  $H^{0,-1}$ are mutually ``complex conjugate"
$\dim(Z)$-dimensional complex vector spaces. This splitting is
$\End^0(Z)$-invariant and the $\End^0(Z)$-module $H^{-1,0}$ is
canonically isomorphic to $\Lie(Z)$. Let
$$\f_H:H_1(Z,\C) \to
H_1(Z,\C)$$ be the $\C$-linear operator in $H_1(Z,\C)$ defined as
follows.
$$\f_H(x) =-x \quad \forall x \in H^{-1,0}; \quad \f_H(x)=0 \quad
\forall x \in H^{0,-1}.$$ Clearly, $\f_H$ commutes with
$\End^0(Z)$ and therefore with $E$. Hence $\f_H$  may be viewed as
an endomorphism of the free $E_{\C}$-module $H_1(Z,\C)$; clearly,
its trace is the tuple
$$(-n_{\sigma})_{\sigma
\in\Sigma_E}\in\prod_{\sigma\in\Sigma_E}\C_{\sigma}=E_{\C}.$$
Suppose $MT=MT_Z \subset \GL_{\Q}(H_1(Z,\Q)))$ is the Mumford-Tate
group of (the rational Hodge structure $H_1(Z,\Q)$ and of) $Z$
(\cite{Deligne,Ribet3,ZarhinIzv}). It is a connected reductive
algebraic $\Q$-group that contains scalars and could be described
as follows (\cite[section 6.3]{ZarhinIzv}). Let $mt\subset
\End_{\Q}(H_1(Z,\Q)))$ be the $\Q$-Lie algebra of $MT$; it is a
reductive algebraic linear $\Q$-Lie algebra  which contains
scalars and  its natural faithful representation in $H_1(Z,\Q)$ is
completely reducible. In addition, $mt$ is the {\sl smallest}
$\Q$-Lie subalgebra in $\End_{\Q}(H_1(Z,\Q)))$  enjoying the
following property: its complexification

$$mt_{\C}=mt\otimes_{\Q}\C \subset \End_{\C}(H_1(Z,\C))$$
contains scalars and $\f_H$. It is well-known that the centralizer
of $MT$ (and therefore of $mt$) in $\End_{\Q}(H_1(Z,\Q))$
coincides with $\End^0(Z)$. This implies that the center $\cc$ of
$mt$ lies in $\CC_Z$. Since $mt$ is reductive, it splits into a
direct sum
$$mt=mt^{ss} \oplus \cc$$
of $\cc$ and a semisimple $\Q$-Lie algebra $mt^{ss}$. Clearly,
$mt$ lies in   $\End_{E}(H_1(Z,\Q))$.

Since $mt^{ss}$ is semisimple and the trace map $\Tr_E$ is a Lie
algebra homomorphism, $\Tr_{E}(mt^{ss})=\{0\}$. Since
$\cc\subset\CC_Z\subset E$, we have $\Tr_{E}(\cc)\subset \CC_Z$
and therefore
$$\Tr_{E}(mt) \subset \CC_Z.$$
This implies easily that
$$\Tr_{E_{\C}}(mt_{\C}) \subset \CC_Z\otimes_{\Q}\C.$$
In particular, since $\f_H \in mt_{\C}$, we have
$\Tr_{E_{\C}}(\f_H) \in \CC_Z\otimes_{\Q}\C$. But
$\Tr_{E_{\C}}(\f_H)=(-n_{\sigma})_{\sigma \in\Sigma}$. This
implies easily that $$(n_{\sigma})_{\sigma
\in\Sigma_E}=-\Tr_{E_{\C}}(\f_H)\in \CC_Z\otimes_{\Q}\C.$$ In
order to prove the second assertion of the theorem, notice that
its assumptions imply that $E/\CC_Z$ is a nontrivial Galois
extension. If $\kappa: E \to E$ is a non-identity element of the
Galois group $\Gal(E/\CC_Z)$ then one may easily check that
$$\CC_{Z}\otimes_{\Q}\C\subset \{(u)_ {\sigma\in\Sigma_E}\in
\prod_{\sigma\in\Sigma_E}\C=E_{\C}\mid
u_{\sigma}=u_{\sigma\kappa} \quad \forall \sigma\}.$$

\section{Cyclic covers and jacobians}
\label{Prelim}

 Throughout this paper we fix an odd prime $p$ and assume
that $K$ is a field of characteristic zero. We fix an algebraic
closure $K_a$ and write $\Gal(K)$ for the absolute Galois group
$\Aut(K_a/K)$. We also fix in $K_a$ a primitive $p$th root of
unity $\zeta$.

Let $f(x) \in K[x]$ be a separable polynomial of degree $n \ge 4$.
We write $\R_f$ for the set of its roots and denote by
$L=L_f=K(\R_f)\subset K_a$ the corresponding splitting field. As
usual, the Galois group $\Gal(L/K)$ is called the Galois group of
$f$ and denoted by $\Gal(f)$. Clearly, $\Gal(f)$ permutes elements
of $\R_f$ and the natural map of $\Gal(f)$ into the group
$\Perm(\R_f)$ of all permutations of $\R_f$ is an embedding. We
will identify $\Gal(f)$ with its image and consider it as a
permutation group of $\R_f$. Clearly, $\Gal(f)$ is transitive if
and only if $f$ is irreducible in $K[x]$.

We refer the reader to \cite{Mortimer,ZarhinCrelle,Klemm,MR} for
the definition and properties of the {\sl heart}
$(\F_p^{\R_f})^{00}$ over the field $\F_p$ of the group $\Gal(f)$
acting on the set $\R_f$. Here we just recall that
$(\F_p^{\R_f})^{00}$ is a finite-dimensional $\F_p$-vector space
provided with a natural structure of $\Gal(f)$-module.

 Let $C=C_{f,p}$ be the smooth projective model of the smooth
affine $K$-curve
            $$y^p=f(x).$$
So $C$ is a smooth projective curve defined over $K$. The
rational function $x \in K(C)$ defines a finite cover  $\pi:C \to
\P^1$ of degree $p$. Let $B'\subset C(K_a)$ be the set of
ramification points.  Clearly, the restriction of $\pi$ to $B'$ is
an {\sl injective} map $B' \hookrightarrow \P^1(K_a)$, whose image
is the disjoint union of $\infty$ and  $\R_f$ if $p$ does {\sl
not} divide $\deg(f)$ and just $\R_f$ if it does. We write
$$B=\pi^{-1}(\R_f)=\{(\alpha,0)\mid \alpha \in \R_f\} \subset B' \subset C(K_a).$$
Clearly, $\pi$ is ramified at each point of $B$ with ramification
index $p$. We have $B'=B$ if and only if $n$ is  divisible by $p$.
If $n$ is not divisible by $p$ then $B'$ is the disjoint union of
$B$ and a single point $\infty':=\pi^{-1}(\infty)$. In addition,
the ramification index of $\pi$ at $\pi^{-1}(\infty)$ is also $p$.
Using Hurwitz's formula, one may easily compute the genus
$g=g(C)=g(C_{p,f})$ of $C$ (\cite[pp.~401--402]{Koo},
\cite[proposition 1 on p. 3359]{Towse}, \cite[p. 148]{Poonen}).
Namely, $g$ is $(p-1)(n-1)/2$ if $p$ does {\sl not} divide $n$ and
$(p-1)(n-2)/2$ if it does.

\begin{rem}
\label{genre}
 Assume that $p$ does not divide $n$ and consider  the plane triangle
(Newton polygon)
$$\Delta_{n,p}:=\{(j,i)\mid 0\le j,\quad 0\le i, \quad pj+ni\le np\}$$
with the vertices $(0,0)$, $(0,p)$ and $(n,0)$. Let $L_{n,p}$ be
the set of integer points in the interior of $\Delta_{n,p}$.
 One may easily check that $g$ coincides with the number
of elements of $L_{n,p}$.  It is also clear that for each
$(j,i)\in L_{n,p}$
$$1\le j \le n-1; \quad 1 \le i \le p-1;\quad p(j-1)+(j+1)\le n(p-i).$$
Elementary calculations (\cite[theorem 3 on p. 403]{Koo}) show
that
$$\omega_{j,i}:=x^{j-1}dx/y^{p-i}=x^{j-1}y^idx/y^p=x^{j-1}y^{i-1} dx/y^{p-1}$$ is a differential of the first
kind on $C$ for each $(j,i) \in L_{n,p}$. This implies easily that
the collection $\{\omega_{j,i}\}_{(j,i)\in L_{n,p}}$ is a basis in
the space of differentials of the first kind on $C$.
\end{rem}

 There is a non-trivial birational
 $K_a$-automorphism of $C$
 $$\delta_p:(x,y) \mapsto (x, \zeta y).$$
Clearly, $\delta_p^p$ is the identity map and
 the set of fixed points of $\delta_p$ coincides with $B'$.

Let $J^{(f,p)}=J(C)=J(C_{f,p})$ be the jacobian of $C$. It is a $g$-dimensional abelian variety defined
 over $K$ and one may view (via Albanese functoriality) $\delta_p$ as an element of
 $$\Aut(C) \subset\Aut(J(C)) \subset \End(J(C))$$
such that
  $\delta_p \ne \I$ but $\delta_p^p=\I$
where $\I$ is the identity endomorphism of $J(C)$. Here $\Aut(C)$
stands for the group of $K_a$-automorphisms of $C$, $\Aut(J(C))$
stands for the group of $K_a$-automorphisms of $J(C)$ and
$\End(J(C))$ stands for the ring of all $K_a$-endomorphisms of
$J(C)$.
 As usual, we write $\End^0(J(C))=\End^0(J^{(f,p)})$ for
the corresponding $\Q$-algebra $\End(J(C))\otimes \Q$.

\begin{lem}
\label{cycl}
 $\I+\delta_p+ \cdots + \delta_p^{p-1}=0$ in $\End(J(C))$.
Therefore the subring $\Z[\delta_p] \subset \End(J(C))$ is
isomorphic to the ring $\Z[\zeta_p]$ of integers in the $p$th
cyclotomic field $\Q(\zeta_p)$. The $\Q$-subalgebra
$\Q[\delta_p]\subset\End^0(J(C))=\End^0(J^{(f,p)})$ is isomorphic
to $\Q(\zeta_p)$.
\end{lem}

{\sl Proof.} See \cite[p.~149]{Poonen},  \cite[p.~458]{SPoonen}.

\begin{rem}
If $K$ contains $\zeta$ then the Galois modules
$(\F_p^{\R_f})^{00}$ and $\ker(\I-\delta_p)$ are canonically
isomorphic (\cite[proposition~6.2]{Poonen},
 \cite[proposition~3.2]{SPoonen}).
\end{rem}

\begin{rem}
\label{nondiv} Recall that $p$ is odd and assume that $n=\deg(f)$
is divisible by $p$ say, $n=pm$ for some positive integer $m$.
Since $n \ge 4$, we conclude that $n \ge 5$.

 Let $\alpha \in K_a$ be a root of $f$ and $K_1=K(\alpha)$ be
the corresponding subfield of $K_a$. We have $f(x)=(x-\alpha)f_1(x)$ with $f_1(x) \in K_1[x]$. Clearly, $f_1(x)$ is a separable polynomial over $K_1$ of   degree $pm-1=n-1 \ge 4$. It is also clear that the polynomials
$$h(x)=f_1(x+\alpha), h_1(x)=x^{n-1}h(1/x) \in K_1[x]$$
are separable of the same degree $pm-1=n-1\ge 4$. The standard
substitution $$x_1=1/(x-\alpha), y_1=y/(x-\alpha)^m$$ establishes
a birational isomorphism between $C_{f,p}$ and a curve
$$C_{h_1}: y_1^p=h_1(x_1)$$ (see \cite[p.~3359]{Towse}). But
$\deg(h_1)=pm-1$ is {\sl not} divisible by $p$. Clearly, this
isomorphism commutes with the actions of $\delta_p$.
\end{rem}

\begin{thm}
\label{handyp3} Suppose $n \ge 4$. Assume that $\Q[\delta_p]$ is a
maximal commutative subalgebra in $\End^0(J^{(f,p)})$. Then the
center $\CC$ of $\End^0(J^{(f,p)})$ is a CM-subfield of
$\Q[\delta_p]$.
\end{thm}

{\sl Proof}. This is theorem 3.8 of \cite{ZarhinCrelle}.

\begin{thm}
\label{bigCM} Suppose $n \ge 4$. Assume that $\Q[\delta_p]$ is a
maximal commutative subalgebra in $\End^0(J^{(f,p)})$. Then
$\End^0(J^{(f,p)})=\Q[\delta_p] \cong\Q(\zeta_p)$ and therefore
 $\End(J^{(f,p)})=\Z[\delta_p] \cong \Z[\zeta_p]$.
\end{thm}

{\sl Proof}. Let  $\CC=\CC_{J^{(f,p)}}$ be the center of
$\End^0(J^{(f,p)})$. We know that $\CC$ is a CM-subfield of
$E:=\Q[\delta_p]$.

Replacing, if necessary, $K$ by  its subfield (finitely) generated
over $\Q$ by all the coefficients of $f$, we may assume that $K$
(and therefore $K_a$) is isomorphic to a subfield of the field
$\C$ of complex numbers. So, $K \subset K_a \subset \C$. We may
also assume that $\zeta=\zeta_p$ and consider $C_{(f,p)}$ as
complex projective curve and its jacobian $J^{(f,p)}$ as complex
abelian variety.

Let $\Sigma=\Sigma_E$ be the set of all field embeddings
$\sigma:E=\Q[\delta_p]\hookrightarrow \C$. We  are going to apply
Theorem \ref{mult}  to $Z=J^{(f,p)}$ and $E=\Q[\delta_p]$. In
order to do that we need to get some information about the
multiplicities
$$n_{\sigma}=\n_{\sigma}(Z,E)=n_{\sigma}(J^{(f,p)},\Q[\delta_p]).$$
Remark \ref{dual} allows us to do it, using the action of
$\Q[\delta_p]$  on the space $\Omega^1(J^{(f,p)})$ of
differentials of the first kind on $J^{(f,p)}$.

Recall that if ${\sigma}':\Q[\delta_p]\hookrightarrow \C$ is the
embedding complex conjugate to $\sigma$ then, by Remark
\ref{conjugate},
$$n_{\sigma}+n_{{\sigma}'}=\frac{2\dim(J^{(f,p)})}{p-1},$$ since
$[\Q[\delta_p]:\Q]=p-1$. Notice also that for each $\sigma:
\Q[\delta_p]\hookrightarrow \C$
$$\Omega^1(J^{(f,p)})_{\sigma}=\{\omega\in \Omega^1(J^{(f,p)})\mid
\delta_p(\omega)=\sigma(\delta_p)\omega\}.$$

In other words, $\Omega^1(J^{(f,p)})_{\sigma}$ is the
eigenspace corresponding to the eigenvalue $\sigma(\delta_p)$
of $\delta_p$ and $n_{\sigma}$ is the multiplicity of the
eigenvalue $\sigma(\delta_p)$.

 Let $i<p$ be a positive integer and
$\sigma_i:\Q[\delta_p]\hookrightarrow \C$ be the embedding which
sends $\delta_p$ to $\zeta^{-i}$. Obviously, the complex conjugate
of $\sigma_i$ coincides with $\sigma_{p-i}$. In addition, for each
$\sigma$ there exists precisely one $i$ such that
$\sigma=\sigma_i$. Clearly, $\Omega^1(J^{(f,p)})_{\sigma_i}$ is
the eigenspace of $\Omega^1(J^{(f,p)})$ attached to the
eigenvalue $\zeta^{-i}$ of $\delta_p$. Therefore $n_{\sigma_i}$
coincides with the multiplicity of the eigenvalue $\zeta^{-i}$.

Let $P_0$ be one of the $\delta_p$-invariant points (i.e., a
ramification point for $\pi$) of $C_{f,p}(K_a)\subset
C_{f,p}(\C)$. Then $$\tau: C_{f,p} \to J^{(f,p)}, \quad P\mapsto
\cl((P)-(P_0))$$ is an embedding of complex algebraic varieties
and it is well-known that the induced map $$\tau^*:
\Omega^1(J^{(f,p)}) \to \Omega^1(C_{f,p})$$ is a $\C$-linear
isomorphism obviously commuting with the actions of $\delta_p$.
(Here $\cl$ stands for the linear equivalence class.) This implies
that $n_{\sigma_i}$ coincides with the dimension of the
eigenspace of $\Omega^1(C_{(f,p)})$ attached to the eigenvalue
$\zeta^{-i}$ of $\delta_p$.

\begin{rem}
\label{eigen}
 Clearly, if for some positive integer $j$ the differential $x^{j-1} dx/y^{p-i}$ lies in
$\Omega^1(C_{(f,p)})$ then it is an eigenvector of $\delta_p$ with
eigenvalue $\zeta^{i}$. Now assume that $p$ does {\sl not} divide
$n$. It follows from Remark \ref{genre} that  each
$n_{\sigma}=n_{\sigma_i}$ could be visualized as the number of
interior integer points in $\Delta_{n,p}$ along the corresponding
(to $p-i$) horizontal line. Elementary calculations show that this
number is $[\frac{ni}{p}]$. This implies that
$n_{\sigma_i}=[\frac{ni}{p}]$ for $1\le i \le p-1$. Then
$n_{\sigma_i}=0$ if and only if $1\le i \le [\frac{p}{n}]$.

Assume, in addition, that $p<n$. Clearly, in this case the
function $i \mapsto n_{\sigma_i}=[\frac{ni}{p}]$ is strictly
increasing.
\end{rem}

\begin{rem}
\label{eigend} Assume that $p$ divides $n$. Then $n \ge 5$ and
$n-1 \ge 4$. Clearly,  $p$ does {\sl not} divide $n-1$. Applying
 Remark \ref{nondiv}, we get a curve $C_{h_1,p}: y_1^p=h_1(x_1)$
 with separable
polynomial $h_1(x_1)$ of degree $n-1$ and  a
$\delta_p$-equivariant birational isomorphism between $C_{f,p}$
and $C_{h_1,p}$. This gives us a $\delta_p$-equivariant
isomorphism
$$\Omega^1(C_{f,p}) \cong \Omega^1(C_{h_1,p}).$$
Applying Remark \ref{eigen} to $C_{h_1,p}$ and $n-1$ (instead of
$C_{f,p}$ and $n$), we conclude that
$n_{\sigma_i}=[\frac{(n-1)i}{p}]$ for $1\le i \le p-1$. Then
$n_{\sigma_i}=0$ if and only if $1\le i \le [\frac{p}{n-1}]$.

Assume, in addition, that $n \ne p$. Clearly, in this case $n-1>p$
and the function $i \mapsto n_{\sigma_i}=[\frac{(n-1)i}{p}]$ is
strictly increasing.
\end{rem}

\begin{prop}
\label{bigp}
\begin{enumerate}
\item[(i)]
let us assume that $p > n$. Then $n_{\sigma}=0$ if and only if
$\sigma=\sigma_i$ with $1 \le i \le [\frac{(p-1)}{n}]$.
\item[(ii)]
let us assume that $p=n$. Then $n_{\sigma}=0$ if and only if
$\sigma=\sigma_i$ with $i=1$.
\end{enumerate}
\end{prop}

{\sl Proof of Proposition \ref{bigp}}. First, assume that $p>n$.
 Clearly, $p$ does {\sl not} divide
$n$ and  therefore $[\frac{p}{n}]=[\frac{p-1}{n}]$. Now the
assertion (i) follows from Remark \ref{eigen}.

 Now assume that $n=p$. By Remark \ref{eigend},
$n_{\sigma_i}=0$ if and only if $1\le i \le [\frac{p}{n-1}]$. But
 $[\frac{p}{n-1}]=[\frac{p}{p-1}]=1$.

\begin{prop}
\label{smallp} Let us assume that $p<n$. If $\sigma, \iota$ are
two embeddings $\Q[\delta_p]\hookrightarrow \C$ then
$n_{\sigma}=n_{\iota}$ if and only if $\sigma=\iota$.
\end{prop}

{\sl Proof of Proposition \ref{smallp}}. First assume that $p$
does {\sl not} divide $n$. Then the assertion follows from (the
last sentence of) Remark \ref{eigen}.

 Now assume that $p$ divides $n$. Then the assertion follows from (the
last sentence of) Remark \ref{eigend}.

\smallskip

{\bf End of the proof of Theorem \ref{bigCM}}. If
$\CC=\Q[\delta_p]$ then we are done, since $\Q[\delta_p]$ is a
maximal commutative subalgebra in $\End^0(J^{(f,p)})$. Assume that
$\CC \ne \Q[\delta_p]$. Our goal is to get a contradiction.

Clearly, $\Q[\delta_p]/\Q$ is a Galois extension. It follows from
Theorem \ref{mult} and Remark \ref{dual} (applied to $Z=
J^{(f,p)}$ and $E=\Q[\delta_p]$) that there exists a non-trivial
field automorphism
 $\kappa: \Q[\delta_p] \to \Q[\delta_p]$ such that
 for all $\sigma\in \Sigma$
$$n_{\sigma}=n_{\sigma\kappa}.$$
 Clearly, there exists an integer $m$ such that
$1<m<p$ and $\kappa(\delta_p)=\delta_p^m$.

First, assume that $n>p$. It follows from Proposition
\ref{smallp} that $\sigma\kappa=\sigma$ which could not be the
case, since $\kappa$ is {\sl not} the identity map. This
contradiction proves the Theorem in the case of $n>p$.

Second, assume that $n=p$. It follows from Proposition
\ref{bigp}(ii) that $\sigma_1\kappa=\sigma_1$ which, by the same
token, leads to a contradiction.

Third, assume that $p>n$. It follows from Proposition
\ref{bigp}(i) that the map $\sigma \mapsto \sigma \kappa$ permutes
the set $\{\sigma_i\mid 1 \le i \le [\frac{(p-1)}{n}]\}$. Since
$\kappa(\delta_p)=\delta_p^m$,
$\sigma_i\kappa(\delta_p)=\zeta^{-im}$. This implies that
multiplication by $m$ in $\F_p^{*}$ leaves invariant the subset
$A:=\{i\bmod p\in \F_p\mid 1  \le i \le [\frac{(p-1)}{n}]\}$. This
implies that $$m=m\cdot 1\le \left[\frac{(p-1)}{n}\right] \le
\frac{(p-1)}{4}.$$ Let us consider the arithmetic progression
consisting  of the $m$ integers $[\frac{(p-1)}{n}]+1, \ldots ,
[\frac{(p-1)}{n}]+m$ with difference $1$. All its elements lie
between $[\frac{(p-1)}{n}]+1$ and
$$\left[\frac{(p-1)}{n}\right]+m  \le 2\left[\frac{(p-1)}{n}\right]\le
2\frac{(p-1)}{4}=\frac{(p-1)}{2}<p-1.$$ Clearly, there exists a
positive integer $r \le m$ such that $[\frac{(p-1)}{n}]+r$ is
divisible by $m$, i.e., there is a positive integer $d$ such that
$md=[\frac{(p-1)}{n}]+r$. Since $[\frac{(p-1)}{n}]\ge m \ge 2$, we
have $d \le [\frac{(p-1)}{n}]$ but $md=[\frac{(p-1)}{n}]+r\le
[\frac{(p-1)}{n}]+m <p-1$. This implies that $A$ is {\sl not}
invariant under multiplication by $m$ which gives the desired
contradiction.

\section{Jacobians and their endomorphism rings}
\label{prf}
 Recall that $K$ is a field of characteristic zero, $K_a$ is
its algebraic closure.
Suppose  $f(x)\in K[x]$ is a polynomial of degree $n \ge 5$
without multiple roots, $\R_f \subset K_a$ is the set of its
roots, $K(\R_f)$ is its splitting field. Let us put
$$\Gal(f)=\Gal(K(\R_f)/K)\subset\Perm(\R_f).$$

\begin{thm}[corollary 5.3  of \cite{ZarhinCrelle}]
\label{handyV} Let $p$ be an odd prime. If $f(x)\in K[x]$ is an
irreducible polynomial of degree $n\ge 5$ and $\Gal(f)=\Sn$ or
$\An$ then $\Q[\delta_p]$ is a maximal commutative subalgebra in
$\End^0(J^{(f,p)})$ and the center of $\End^0(J^{(f,p)})$ is a
CM-subfield of $\Q[\delta_p]$.
\end{thm}

Combining Theorems \ref{handyV} and \ref{bigCM}, we obtain the
following statement.

\begin{thm}
\label{final} Let $p$ be an odd prime. If $f(x)\in K[x]$ is an
irreducible polynomial of degree $n\ge 5$ and $\Gal(f)=\Sn$ or
$\An$ then
 $\End^0(J^{(f,p)})=\Q[\delta_p]$ and therefore
$\End(J^{(f,p)})=\Z[\delta_p]\cong \Z[\zeta_p]$.
\end{thm}

 Clearly, Theorem \ref{endo}  is a special case of Theorem
 \ref{final}.

\begin{ex}
Suppose $L=\C(z_1, \cdots , z_n)$ is the field of rational
functions in $n$ independent variables $z_1, \cdots , z_n$ with
constant field $\C$ and $K=L^{\Sn}$ is the subfield of symmetric
functions. Then $K_a=L_a$ and $$f(x)=\prod_{i=1}^n(x-z_i) \in
K[x]$$ is an irreducible polynomial over $K$ with Galois group
$\Sn$. Let $C$ be a smooth projective model of the $K$-curve
$y^p=f(x)$ and $J(C)$ its jacobian. It follows from  Theorem
\ref{final} that if $n\ge 5$ then the ring of $L_a$-endomorphisms
of $J(C)$ is $\Z[\zeta_p]$. In particular, the abelian variety
$J(C)$ is absolutely simple. When $p=3$ and $3\mid n$ the absolute
simplicity of $J(C)$ was proven in (\cite[p.~107]{V}).
\end{ex}

\begin{ex}
Let $h(x) \in \C[x]$ be a  {\sl Morse polynomial} of degree $n \ge
5$. This means that the derivative $h'(x)$ of $h(x)$ has $n-1$
distinct roots $\beta_1, \cdots \beta_{n-1}$  and $h(\beta_i) \ne
h(\beta_j)$ while $i\ne j$. (For example, $x^n-x$ is a Morse
polynomial.) Let $K=\C(z)$ be the field of rational functions in
variable $z$ with constant field $\C$ and $K_a$ its algebraic
closure. Then a theorem of Hilbert  (\cite[theorem~4.4.5,
p.~41]{Serre}) asserts that the Galois group of $h(x)-z$ over
$k(z)$ is $\Sn$. Let $C$ be a smooth projective model of the
$K$-curve $y^p=h(x)-z$ and $J(C)$ its jacobian. It follows from
Theorem \ref{final} that the ring of $K_a$-endomorphisms of $J(C)$
is
 $\Z[\zeta_p]$. In particular, the abelian variety $J(C)$ is absolutely simple.
\end{ex}

We refer the reader to
\cite{ZarhinM,ZarhinMRL2,ZarhinCrelle,ZarhinMMJ}
 for the definition and basic properties of {\sl very
simple} representations.

\begin{thm}
\label{handysup} Suppose $p$ is an odd prime, $n\ge 5$ and  $K$
contains a primitive $p$th root of unity. If the $\Gal(f)$-module
$(\F_p^{\R_f})^{00}$ is very simple then $\Q[\delta_p]$ coincides
with its own centralizer in $\End^0(J^{(f,p)})$.
\end{thm}

{\sl Proof}. See theorem 5.2 of \cite{ZarhinCrelle}.

\begin{thm}
Suppose $p$ is an odd prime, $n\ge 5$  and  $K$ contains a
primitive $p$th root of unity. If the $\Gal(f)$-module
$(\F_p^{\R_f})^{00}$ is very simple then
$\End^0(J^{(f,p)})=\Q[\delta_p]$ and therefore
$\End(J^{(f,p)})=\Z[\delta_p]\cong \Z[\zeta_p]$.
\end{thm}

{\sl Proof.} It is an immediate corollary of Theorem
\ref{handysup} combined with Theorem \ref{bigCM}.

\section{Corrigendum to \cite{ZarhinCrelle}}
\label{Crelle}
 Remark 2.1 on p. 94, the last assertion. In
general, it is not necessarily true that $G$ is doubly transitive
\cite[Beispiel~2c]{Klemm}, \cite{MR}. However, it becomes true if
one assumes additionally that either $p$ does not divide $n$ or
$G$ is transitive and $p$ is an odd number dividing $n$
(\cite[Satz~4a and Satz~11]{Klemm},
\cite[lemma~2.4]{ZarhinCrelle}).

Lemma 2.4  on p. 95. Its assertion is essentially contained in
Satz 4a of \cite{Klemm}.

Remark 2.5 on p. 95. Its assertion is essentially  Hilffsatz 3b of
\cite{Klemm}.

Sections 1, 3 and 5. Everywhere $\Q(\delta_p)$ means
$\Q[\delta_p]$. (However, it does not make a difference, since
$\Q[\delta_p]$ is a field.)


\begin{thebibliography}{99}
\bibitem{Deligne}P.~Deligne, Hodge cycles on abelian varieties (notes by
J.\/S.~Milne). Lecture Notes in Math., vol. {\bf 900}
(Springer-Verlag, 1982), pp. 9--100.

\bibitem{JN} J. de Jong and R. Noot, Jacobians with complex
multiplications. In: {\em Arithmetic algebraic geometry} (eds. G.
van der Geer, F. Oort and J. Steenbrink). Progress in Math., vol.
{\bf 89} (Birkh\"auser, 1991), pp. 177--192.

\bibitem{Klemm} M. Klemm, \"Uber die Reduktion von
Permutationsmoduln. {\em Math. Z.} {\bf 143} (1975), 113--117.



\bibitem{Koo} J. K. Koo, On holomorphic differentials of some algebraic function field of one variable over
 C. {\em Bull. Austral. Math. Soc.} {\bf 43} (1991), 399--405.

\bibitem{MZ} B. Moonen and Yu. G. Zarhin, Weil classes on abelian varieties. {\em J. reine angew. Math.} {\bf 496} (1998), 83--92.

\bibitem{Mortimer} B. Mortimer, The modular permutation
representations of the known doubly transitive groups. {\em Proc.
London Math. Soc.} (3) {\bf 41} (1980), 1--20.


\bibitem{MumfordAV} D. Mumford, {\em Abelian varieties}, 2nd edn (Oxford University Press, 1974).

\bibitem{MR} P. M. Neumann, Review of \cite{Klemm}, {\em MR} 52 \#544.

\bibitem{Ribet3} K. Ribet, Hodge classes on certain abelian varieties. {\em Amer. J. Math.} {\bf 105} (1983), 523--538.


\bibitem{Poonen} B. Poonen and E. Schaefer, Explicit descent for Jacobians of cyclic covers of
 the projective line. {\em J. reine angew. Math.} {\bf 488} (1997), 141--188.

\bibitem{SPoonen} E. Schaefer, Computing a Selmer group of a Jacobian using functions on the curve.
{\em Math. Ann.} {\bf 310} (1998), 447--471.

\bibitem{Schur} I. Schur, Gleichungen ohne Affect. {\em Sitz. Preuss. Akad. Wiss.} 1930,
 {\em Physik-Math. Klasse} 443--449 (=Ges. Abh. III, 191--197).

\bibitem{Serre} J.-P. Serre, {\em Topics in Galois Theory} (Jones and Bartlett Publishers, 1992).


\bibitem{Towse} C. Towse, Weierstrass points on cyclic covers of the projective line. {\em Trans. Amer. Math. Soc.}
 {\bf 348} (1996), 3355-3377.

\bibitem{V} H. V\"olklein, Cyclic covers of $\P^1$  and
Galois actions on their division points. {\em Contemporary Math.}
 {\bf 186} (1994), 91--107.

\bibitem{ZarhinIzv} Yu. G. Zarhin, Weights of simple Lie algebras in the cohomology of
algebraic varieties. {\em Izv. Akad. Nauk SSSR Ser. Mat.} {\bf 48}
(1984), 264--304; English translation: {\em Math. USSR Izv.} {\bf
24} (1985), 245 - 281.

\bibitem{Zarhin} Yu. G. Zarhin, Hyperelliptic jacobians without complex multiplication.
 {\em Math. Res. Letters} {\bf 7} (2000), 123--132.

\bibitem{ZarhinM} Yu. G. Zarhin, Hyperelliptic jacobians and modular
representations. In: {\em Moduli of abelian varieties} (eds. C.
Faber, G. van der Geer and F. Oort). Progress in Math., vol. {\bf
195} (Birkh\"auser, 2001), pp. 473--490.

\bibitem{ZarhinMRL2}  Yu. G. Zarhin, Hyperelliptic jacobians without
complex multiplication in positive characteristic. {\em Math. Res.
Letters} {\bf 8} (2001), 429--435.

\bibitem{ZarhinCrelle} Yu. G. Zarhin, Cyclic covers of the projective line, their jacobians and endomorphisms.
 {\em J. reine angew. Math.} {\bf 544} (2002), 91--110.

\bibitem{ZarhinSb} Yu. G. Zarhin, Endomorphism rings of
certain jacobians in finite characteristic. {\em Matem. Sbornik}
{\bf 193} (2002), issue 8 (Russian), to appear.

\bibitem{ZarhinMMJ} Yu. G. Zarhin, Very simple $2$-adic representations
 and hyperelliptic jacobians; available at http://arXiv.org/abs/math.AG/0109014; {\em Moscow Math. J.}
  {\bf 2} (2002), issue 2, to appear.

\end{thebibliography}
\end{document}